\documentclass[12pt]{amsart}
\usepackage{amsmath,amscd,amssymb,amsfonts}
\usepackage[all]{xy}
\usepackage{lineno}

\newcommand{\ver}{\today}
\newcommand{\bC}{{\mathbb C}}

\newcommand{\bN}{{\mathbb N}}
\newcommand{\bP}{{\mathbb P}}

\newcommand{\bQ}{{\mathbb Q}}
\newcommand{\bZ}{{\mathbb Z}}

\newcommand{\cF}{{\mathcal F}}

\newcommand{\Z}{\mathbb{Z}}
\newcommand{\M}{M_{f,0}}

\newcommand{\Q}{\mathbb{Q}}

\newcommand{\C}{\mathbb{C}}

\newcommand{\Sf}{\mathcal{O}}
\newcommand{\X}{\chi}

\theoremstyle{plain}
\newtheorem{thm}{Theorem}[section]
\newtheorem{cor}[thm]{Corollary}

\newtheorem{prop}[thm]{Proposition}

\theoremstyle{definition}

\newtheorem{example}[thm]{Example}

\title[Spectrum of plane curve arrangements]{Spectrum of projective plane curve arrangements with ordinary singularities}
\author{Youngho Yoon}
\address{Department of Mathematical Sciences, Seoul National University GwanAkRo 1, Gwanak-Gu, Seoul 08826, Korea} 
\email{nsyyh@snu.ac.kr}

\date{\ver}


\begin{document}

\begin{abstract}
The Hodge spectrum is an important analytic invariant of singularities encoding the Hodge filtration and the monodromy of the Milnor fiber. However, explicit formulas exist in only a few cases. In this article the main result is a combinatorial formula for the Hodge spectrum of any homogeneous polynomials in three variables whose zero locus is a projective curve arrangement having only ordinary multiple points. 
\end{abstract}

\maketitle

\section{Introduction}
The spectrum $Sp(f)$ of the germ of a hypersurface singularity $f:(\bC^n,\ 0)\rightarrow (\bC,\ 0)$ is a fractional Laurent polynomial $$Sp(f)=\sum_{\alpha\in\bQ} n_{f,\alpha} t^\alpha$$ with $n_{f,\alpha}\in \bZ$. The formal definition of $Sp(f)$ involves the Hodge filtration and the monodromy on the cohomology of the Milnor fiber $M_{f,0}$ of $f$ (see Section \ref{multiplicity}).

The spectra of reduced hyperplane arrangements are known to be combinatorial. Combinatorial formulas for the spectrum are given in \cite{Bu-S} and \cite{YY} in the cases of $n=3$ and $n=4$. One may ask similarly for the spectra of homogeneous polynomials defining projective curve arrangements having only ordinary multiple points (i.e. locally biholomorphic to a line arrangement). We answer this question by giving a combinatorial formula. The formula is calculated in a general setting, where the arrangements might not be reduced arrangements. We also generate formulas for some special cases. 

In this paper, we will follow the convention $\binom{t}{k}=t(t-1)\cdots(t-k+1)/k!$ for $k\in \bN$ and  any $t$. Also, $\lfloor \beta \rfloor:=\max\{z\in \bZ| z \leq\beta\}$ and $\lceil \beta \rceil:=\min\{z\in \bZ| \beta \leq z\}$ for $\beta \in \bQ$. The notation $\delta_{a,b}$ means $1$ if $a=b$, and $0$ otherwise.

\begin{thm}\label{c3}
Let  $f=\prod_{l\in L} f_l^{m_l}$ be the defining equation of a curve arrangement with irreducible component corresponding to $f_l$ in $\bP^2$. Assume the reduced arrangement defined by $f_{red}:=\prod_{l\in L}{f_l}$ has only ordinary singularities. Let $m_{v,l}$, $m_{v,red}$, and $m_v$ be the multiplicities of $f_l$, $f_{red}$, and $f$ respectively at an ordinary multiple point $v \in V$ of $f_{red}$.  Also, let $d_l$, $d_{red}$, and $d$ be the degrees of $f_l$, $f_{red}$, and $f$ respectively. Set $L_v:=\{ l \in L  |  f_l(v)=0\} \subset L$ for the set of irreducible curves containing $v\in V$.
Then the following formulas hold for $j\in \{1,\cdots,d\}:$
\begin{align*}
n_{f,\frac{j}{d}}=& \binom{d_{red}-u_{0,j} -1}{2}-\sum_{v \in V}\binom{m_{v,red} - u_{v,j} -1}{2},\\
n_{f,1+\frac{j}{d}}&=\left( (u_{0,j}-1)(d_{red}-u_{0,j}-1) +\sum_{l\in L}\binom{d_l}{2} \right)\\
	&-\sum_{v\in V}\left( u_{v,j}(m_{v,red}-u_{v,j}-1)+ \sum_{l\in L_v}\binom{m_{v,l}}{2} \right), \text{ and}&\\
n_{f,2+\frac{j}{d}}&=\binom{u_{0,j}-1}{2}-\sum_{v\in V} \binom{u_{v,j}}{2}-\delta_{j,d}.
\end{align*}
Otherwise $n_{f,\alpha}=0$. Here $u_{0,j} = \sum_{l\in L} \lceil j m_l/d \rceil d_l-j$ and $u_{v,j} = \sum_{l\in L_v} \lceil j m_l/d\rceil m_{v,l} - \lceil j m_v/d\rceil$. In the case $f=f_{red}$ we have the following formulas for $j\in \{1,\cdots,d\}:$
\begin{align*}
n_{f,\frac{j}{d}}=& \binom{j-1}{2}-\sum_{v\in V}\binom{\lceil j m_v/d\rceil -1}{2},\\
n_{f,1+\frac{j}{d}}&= \left((j-1)(d-j-1)+\sum_{l\in L}\binom{d_l}{2}\right)\\
	&-\sum_{v\in V}\left((\lceil j m_v/d\rceil -1)(m_v-\lceil j m_v/d\rceil)+\sum_{l \in L_v}\binom{m_{v,l}}{2}\right),\text{ and}
\end{align*}
\begin{align*}	
n_{f,2+\frac{j}{d}}&=\binom{d-j-1}{2}-\sum_{v \in V}\binom{m_v-\lceil j m_v/d\rceil}{2}-\delta_{j,d}.
\end{align*}
Otherwise $n_{f,\alpha}=0$.

\end{thm}

We remark that the set $V$ of singular points of $f_{red}$ can be replaced by any subset of $V$ containing the set of non simple normal crossing singular points.

Notice that $d_{red}=\sum_l d_l$, $d=\sum_l m_l d_l$, $m_{v,red}=\sum_l m_{v,l}$, and $m_v=\sum_l m_l m_{v,l}$. Also, $m_l=1$, $d=d_{red}$, $m_v=m_{v,red}$, $u_{0,j}=d-j$, and $u_{v,j}=m_v-  \lceil j m_v/d\rceil$ if $f=f_{red}$. Therefore, all the input quantities appearing in Theorem \ref{c3} depends only on the numbers $m_l$, $d_l$, and $m_{v,l}$. In particular $m_{l}$ is not needed if $f$ is reduced. We can encode these information into a decorated directed graph consisting of two sets of vertexes $L$ and $V$ with the directed edges corresponding to inclusions. We decorate the vertexes in $L$ by $(m_l, d_l)$ and edges by $m_{v,l}$. The spectrum solely depends on $m_l$, $d_l$, $m_{v,l}$ and the inclusion relations, which are encoded in this graph. These data is combinatorial in a strong sense.

\begin{cor}\label{comb}
The Hodge spectrum of a homogeneous polynomial $f$ in three variables is determined by combinatorial data if its reduced polynomial defines a projective curve arrangement with only ordinary singularities.  

\end{cor}

In Section \ref{form2} we give combinatorial spectrum formulas in some special cases. Theses special formulas will be useful in practice. We review the definition of the Hodge spectrum and an algorithm calculating the spectrum of homogeneous polynomials in Section \ref{spectrum}. The proof of the main theorem will be given in Section \ref{hyper}

\noindent
{\bf Acknowledgements}: The author would like to thank Nero Budur and Alexander Dimca for suggesting this problem to me, as well as Gabriel C. Drummond-Cole and Xia Liao for useful comments. This work was supported by BK21 PLUS SNU Mathematical Sciences Division and the National Research Foundation of Korea(NRF) grant funded by the Korea government(MSIP) (No. NRF-2017R1C1B1005166).

\section{More formulas} \label{form2}

A necessary and sufficient condition for $\sum_{l\in L} \binom{m_{v,l}}{2} = 0$ in Theorem \ref{c3} is that all the curves defined by $f_l$ are smooth since $m_{v,l}=1$. It is equivalent to say that the homogeneous polynomial $f_l$ has only an isolated singularity at $0 \in \bC^3$ if it is not a linear form.

\begin{cor}\label{s3}
With the same notation as Theorem \ref{c3} we assume  $f_l$ defines a smooth curve in $\bP^2$ for each $l$. 
Then we have the following formulas for $j\in \{1,\cdots,d\}$:
\begin{align*}
n_{f,\frac{j}{d}}=& \binom{d_{red}-u_{0,j} -1}{2}-\sum_{v \in V}\binom{m_{v,red} - u_{v,j} -1}{2},\\
n_{f,1+\frac{j}{d}}&= (u_{0,j}-1)(d_{red}-u_{0,j}-1) +\sum_{l \in L}\binom{d_l}{2} -\sum_{v \in V} u_{v,j} (m_{v,red}-u_{v,j}-1),\\
\text{ and}&\\
n_{f,2+\frac{j}{d}}&=\binom{u_{0,j}-1}{2}-\sum_{v \in V} \binom{u_{v,j}}{2}-\delta_{j,d}.
\end{align*}
Otherwise $n_{f,\alpha}=0$. Here $u_{0,j} = \sum_{l\in L} \lceil j m_l/d \rceil d_l-j$ and $u_{v,j} = \sum_{ l \in L_v } \lceil j m_l/d\rceil - \lceil j m_v/d\rceil$. 
\end{cor}

If $d_l=1$ for all $l$ in Theorem \ref{c3} we get a line arrangement. Also, this condition is equivalent to $\sum_{l \in L} \binom{d_l}{2}=0$. It implies  $\sum_{l \in L} \binom{m_{v,l}}{2} = 0$ since each hyperplane is smooth.

\begin{cor}\label{p3}
With the same notation as Theorem \ref{c3} we assume  $f$ defines a hyperplane arrangement. 
Then we have the following formulas for $j\in \{1,\cdots,d\}$:
\begin{align*}
n_{f,\frac{j}{d}}=& \binom{d_{red}-u_{0,j} -1}{2}-\sum_{v\in V}\binom{m_{v,red} - u_{v,j} -1}{2},\\
n_{f,1+\frac{j}{d}}&= (u_{0,j}-1)(d_{red}-u_{0,j}-1) -\sum_{v\in V} u_{v,j} (m_{v,red}-u_{v,j}-1),\text{ and}\\
n_{f,2+\frac{j}{d}}&=\binom{u_{0,j}-1}{2}-\sum_{v\in V} \binom{u_{v,j}}{2}-\delta_{j,d}.
\end{align*}
Otherwise $n_{f,\alpha}=0$. Here $u_{0,j} = \sum_{l \in L} \lceil j m_l/d \rceil-j$ and $u_{v,j} = \sum_{ l \in L_v } \lceil j m_l/d\rceil - \lceil j m_v/d\rceil$. 
\end{cor}

The formulas for reduced hyperplane arrangements are given in \cite{Bu-S}.

We can consider a single irreducible curve which has only ordinary multiple points. In this case the index set $L=\{ l \}$ of irreducible components has a single element in Theorem \ref{c3}. Put $m_l=m$. Then we have the equalities: $d_l=d_{red}$, $m_{v,l}=m_{v,red}$, $m_v=m m_{v,red}$ and $d=m d_{red}$.

\begin{cor}\label{k1}
Let $f$ be a degree $d$ homogeneous polynomial defining a degree $d$ irreducible curve in $\bP^2$. Assume the curve has only ordinary multiple points. Let $m_{v}$ be the multiplicity of $f$ at each ordinary multiple point $v$ of the curve. The spectrum of $f^m$ for $m\in \bN$ is given by the following formulas for $j\in \{1,\cdots,md\}$:
\begin{align*}
n_{f^m,\frac{j}{md}}=& \binom{d-u_{0,j} -1}{2}-\sum_{v \in V}\binom{m_v - u_{v,j} -1}{2},\\
n_{f^m,1+\frac{j}{md}}&=\left( (u_{0,j}-1)(d-u_{0,j}-1) +\binom{d}{2} \right)-\sum_{v \in V}\left( u_{v,j} (m_v-u_{v,j}-1)+ \binom{m_{v}}{2} \right),\\
\text{ and}&\\
n_{f^m,2+\frac{j}{md}}&=\binom{u_{0,j}-1}{2}-\sum_{v \in V} \binom{u_{v,j}}{2}-\delta_{j,md}.
\end{align*}
Otherwise $n_{f^m,\alpha}=0$. Here $u_{0,j} = \lceil j /d \rceil d-j$ and $u_{v,j} = \lceil j /d\rceil m_{v} - \lceil j m_v/d\rceil$. 

\end{cor}

Applying the Thom-Sebastiani formula we can calculate the Hodge spectrum of functions defined on lower dimensional spaces (see \cite{Kul}-II (8.10.6) and \cite{YY} Section 5). 

\begin{cor}\label{l2}
Let $f=\prod_{l \in L} f_l^{m_l}: (\bC^2,0) \to (\bC, 0)$ be a product of linear forms $f_l$ in two variables with degree $d$. Also, let $d_{red}$ be the number of different linear forms in $f$. Then we have the following formulas for $j\in \{1,\cdots,d\}$,
\begin{align*}
n_{f,\frac{j}{d}}=& d_{red}- \sum_{l \in L} \lceil j m_l/d\rceil +j-1 \text{ and}\\
n_{f,1+\frac{j}{d}}&=\sum_{l \in L} \lceil j m_l/d\rceil -j-1+\delta_{j,d}.\\
\end{align*}
Otherwise $n_{f,\alpha}=0$. 
\end{cor}

We remark that Corollary \ref{k1} and Corollary \ref{l2} cover the cases of homogeneous isolated singularities.  The Hodge spectrum formula of degree $d$ homogeneous isolated singularities is known as  
$$Sp(f)=\left( \frac{t^{1/d}-t}{1-t^{1/d}} \right)^n.$$
We can check the formulas in Corollary \ref{k1} and Corollary \ref{l2} give the same formula in the cases $n=3$ and $n=2$.  

\begin{example}\label{ex}Let $f_1=x$, $f_2=y$, and $f_3=x^2 z+y^2 z -y^3$.  
	\begin{enumerate}
		\item Consider $f:=f_3$. It has a unique singularity at $(0;0;1)\in \mathbb{P}^2$.
		By Theorem \ref{c3} or Corollary \ref{k1} 
		$$Sp(f)=t+2 t^{4/3}+2 t^{5/3}.$$
		\item In the case $f:=f_1 f_3=x^3z+xy^2z+xz^3$ we have two singular point $v_1:=(0;0;1)$ and $v_2:=(0;1;1)$ in $\mathbb{P}^2$. At these points we have $m_{v_1,1}=1$, $m_{v_1,3}=2$, $m_{v_1}=3$, $m_{v_2,1}=1$, $m_{v_2,3}=1$, $m_{v_2}=2$, $d_1=1$, $d_3=3$ and $d=4$. Applying the formula for $f=f_{red}$  in Theorem \ref{c3} we get
		$$Sp(f)=2t+2 t^{5/4}+2 t^{6/4}+2 t^{7/4}-t^2.$$
		\item If $f:=f_1^{m_1}f_2^{m_2}=x^{m_1}y^{m_2}$ its reduced form $f_{red}=xy$ defines a simple normal crossing line arrangement with a singular point $(0;0;1)\in \mathbb{P}^2$. Let $G.C.D.(m_1,m_2)=g$. We can calculate $Sp(f)$ by using Theorem \ref{c3}, Corollary \ref{p3} or Corollary \ref{l2} with Thom-Sebastiani formula. In the case of $g=1$, $Sp(f)=-t^2$. For a general $f=\left( x^{m_1/g}y^{m_2/g}\right)^g$
			\begin{align*}
				Sp(f)=-\frac{t^{1+\frac{1}{g}}(1-t)^2}{1-t^{1/g}}-t^3.
			\end{align*}
	\end{enumerate}
\end{example}

\section{Spectrum}\label{spectrum}

\subsection{Milnor fiber and Hodge spectrum}\label{multiplicity}
Let $f: (\bC^n,0) \rightarrow (\bC,0)$ be the germ of a non-zero holomorphic function. Then the {\it Milnor fiber} at the origin is 
$$M_{f,0}=\{z\in \bC^n \ | \ |z|<\epsilon \; and\; f(z)=\delta\} \text{ for } 0<|\delta|\ll\epsilon\ll 1.$$ 
Refer \cite{Mil} to see the work of John Milnor on this fiber.

The cohomology groups $H^*(M_f,\bC)$ carry canonical mixed Hodge structures such that the semi-simple part $T_s$ of the monodromy acts as an automorphism of finite order of these mixed Hodge structures (see \cite{St}-12.1.3). The eigenvalues $\lambda$ of the monodromy action on $H^*(M_f,\bC)$ are roots of unity. We define the {\it spectrum multiplicity} of $f$ at $\alpha\in \Q$ to be 
$$n_{f,\alpha}=\sum_{j\in \Z} (-1)^{j-n+1}\dim Gr_F^p \tilde{H}^{j}(\M,\C)_\lambda$$
$$\text{with $p=\lfloor n-\alpha \rfloor$, $\lambda=\exp(-2\pi i\alpha ),$}$$
where $\tilde{H}^j(\M,\C)_\lambda$ is the $\lambda$-eigenspace of the reduced cohomology under $T_s$ and $F$ is the Hodge filtration.
The {\it Hodge spectrum} of the germ $f$ is the fractional Laurent polynomial
$$Sp(f):=\sum_{\alpha \in \Q} n_{f,\alpha}t^\alpha.$$

In the case of isolated singularities the {\it Milnor number} $\mu_f$ is defined to be  $\dim\tilde{H}^{n-1}(M_f,\bC)$. This number is the same as $Sp(f)(1)$ by the definition of spectrum and the cohomology vanishing property of isolated singularities. 

\subsection{Spectrum of homogeneous polynomials} \label{homo}
Assume that $f$ is homogeneous of degree $d$. Then we can consider the divisor $Z\subset \bP^{n-1}=:Y$ defined by $f$. Let $\rho : \tilde{Y}\rightarrow Y$ be an embedded resolution of $Z$ inducing an isomorphism over $Y\backslash Z$. We have a divisor $\tilde{Z}:=\rho ^*Z$ with normal crossings on $\tilde{Y}$. Set $\tilde{Z}=\sum_{w \in W} m_w E_w$ where $E_w$ are the irreducible components with multiplicity $m_w$. Let $\tilde{H}$ be the total transform of a general hyperplane $H$ of $Y$. Then the eigenvalues of the monodromy action are $d$-th roots of unity (see \cite{Bu-HH}-4) and we have the following formula for spectrum multiplicities (see \cite{Bu-S}-1.5).

\begin{prop}\label{neuler} For $\alpha =n-p-\frac{i}{d}\in(0,n]$ with $p\in\Z$ and $i\in[0,d-1]\cap \Z$
	\begin{equation*}
	n_{f,\alpha}=(-1)^{p-n+1}\left(  \X\left(\tilde{Y},\Omega_{\tilde{Y}} ^p (\log \tilde{Z})\bigotimes_{\Sf_{\tilde{Y}}} \Sf_{\tilde{Y}} \left( -i\tilde{H}+\sum_{w \in W} \lfloor im_w /d\rfloor E_w \right) \right) -\delta_{\alpha,n}\right).
	\end{equation*} 
\end{prop}
Here we subtract $\delta_{\alpha, n}$ since $n_{f,\alpha}$ is defined using reduced cohomology. 

	Let $U$ be a divisor on $\tilde{Y}$. The class $u:=[U] \in H^2(\tilde{Y}) $ can be written as $u=-u_0 [\tilde{H}]+\sum_{w \in W} u_w [E_w]$ where $u_0, u_w\in \bZ$. Set $\cF_p(U):=\Omega_{\tilde{Y}} ^{p} (\log \tilde{Z})\otimes_{\Sf_{\tilde{Y}}}\Sf_{\tilde{Y}}(U)$ and consider the function $\X_p : H^2(\tilde{Y})\rightarrow \bZ$ defined by $\X_p(u):=\X(\tilde{Y},\cF_p(U))$ for each $p\in\bZ$. By Proposition \ref{neuler}, $\X_p\left(-i [\tilde{H}]+\sum_{w \in W} \lfloor im_w /d\rfloor [E_w] \right)$ gives  $n_{f,n-p-\frac{i}{d}}$. By putting $k:=n-p-1$ and $j:=d-i$,
\begin{equation}\label{euler}
n_{f,k+\frac{j}{d}}=(-1)^{-k}\left( \X_{n-1-k}\left( (d-j) [\tilde{H}]+\sum_{w \in W} \lfloor (d-j)m_w /d\rfloor [E_w] \right) -\delta_{k+\frac{j}{d},n} \right).
\end{equation}
We calculate $\X_p(u)$ by using Hirzebruch-Riemann-Roch,
\begin{equation}\label{Xp}
\X_p(u)=\int_{\tilde{Y}} ch\left(\cF_p (U) \right) \cdot Td(\tilde{Y})=\int_{\tilde{Y}} ch\left(\wedge ^p \left( \Omega_{\tilde{Y}} ^{1} (\log \tilde{Z})\right) \right)\cdot ch(U)  \cdot Td(\tilde{Y}) ,
\end{equation}
 where $ch(\cF_p (U))$ is the Chern character and $Td(\tilde{Y})$ is the Todd class of the tangent bundle $T\tilde{Y}$. Also,  we will denote $c(\tilde{Y}):=c(T\tilde{Y})$ and $ch(\tilde{Y}):=ch(T\tilde{Y})$. 

Chern classes calculate the Chern character and Todd classes. For a given vector bundle $E$ of rank $r$, the following is well known (see \cite{Fu}-3.2),
\begin{equation}\label{generalch}
ch(E)=r+c_1(E)+\frac{1}{2}\left(c_1(E)^2-2c_2(E)\right)+\frac{1}{6}\left(c_1(E)^3 -3c_1(E)c_2(E)+3c_3(E)\right)+\cdots,\\
\end{equation}
and
\begin{equation}\label{td}
Td(E)=1+\frac{1}{2}c_1(E)+\frac{1}{12}\left(c_1(E)^2+c_2(E)\right)+\frac{1}{24}\left(c_1(E)c_2(E)\right)+\cdots.
\end{equation}

Hence, we need to calculate $c(\wedge ^p (\Omega_{\tilde{Y}}^1 (\log \tilde{Z})) )$ to get $ch(\wedge ^p (\Omega_{\tilde{Y}}^1 (\log \tilde{Z})))$. It can be calculated from the following formula. Let $A$ be a rank $r$ vector bundle and write the Chern polynomial $c_t(A):=\sum_{j=0}^r c_j(A)t^r=\prod_{i=1}^r(1+x_i t)$ where $x_1, \cdots, x_r$ are Chern roots of $A$. Then we have the Chern polynomial 
$$c_t(\wedge^p A)=\prod_{1\leq {i_1}<\cdots<{i_p}\leq r}\left(1+(x_{i_1}+\cdots+x_{i_p})t\right).$$

We give an example for our case.
\begin{example}\label{wedge} Let $A$ have rank $2$. \\
	\begin{align*}
		c(\wedge^0 A)&=1,\\
		c(\wedge^1 A)&=1+c_1(A)+c_2(A),\\
		c(\wedge^2 A)&=1+c_1(A)
	\end{align*}
	and
	\begin{align*}
		ch(\wedge^0 A)&=1,\\
		ch(\wedge^1 A)&=2+c_1(A)+\frac{1}{2}\left(c_1(A)^2-2c_2(A)\right),\\
		ch(\wedge^2 A)&=1+c_1(A)+\frac{1}{2}c_1(A)^2
	\end{align*}
\end{example}

We should calculate $c(A)$ for $A=\Omega_{\tilde{Y}}^1(\log\tilde{Z})$. We have the following short exact sequences
\begin{align*}
0\rightarrow \Omega_{\tilde{Y}}^1\rightarrow\Omega_{\tilde{Y}}^1(\log\tilde{Z})\rightarrow \bigoplus_{w \in W} i_{w*} \Sf_{E_w}\rightarrow 0
\end{align*}
 and 
 \begin{align*}
0\rightarrow \Sf_{\tilde{Y}}(-E_W)\rightarrow \Sf_{\tilde{Y}} \rightarrow i_{w*} \Sf_{E_w}\rightarrow 0 
\end{align*}
where $i_w : E_w \to \tilde{Y}$ is the inclusion. This induces
\begin{equation}\label{log}
c(\Omega_{\tilde{Y}}^1(\log\tilde{Z}))=c(\Omega_{\tilde{Y}}^1)\prod_{w \in W}c(i_{w*}\Sf_{E_w})=c(\Omega_{\tilde{Y}}^1)\prod_{w \in W}c(\Sf_{\tilde{Y}}(-E_w))^{-1}.
\end{equation}

Now we apply this algorithm to prove the main theorem.

\section{Proof of Main theorem}\label{hyper}

Let  $f :\bC^3\rightarrow \bC$ be a homogeneous polynomial. We use the same notation and assumptions as in Theorem \ref{c3}. Let $Z=\sum_{l \in L} m_l Z_l \subset \bP^2=:Y$ be the divisor defined by $f$ and $Z_l$ be the irreducible divisor defined by $f_l$. The blowing-up $\rho:\tilde{Y} \to Y$ centered at the set of points $V:=\{v\}$ is an embedded resolution  of  $Z\subset Y$. The pullback $\tilde{Z}:=\rho^*{Z}$ of $Z$ is decomposed as $\tilde{Z}=\sum_{l \in L} m_l E_l +\sum_{v \in V} m_v E_v$ where $E_l$ is the irreducible component corresponding to the strict transform of the irreducible component $Z_l$ and $E_v$ is the exceptional divisor coming from blowing up at $v\in V$. Notice that the set $L \cup V$ is the set $W$ in Section \ref{homo}.

The cohomology ring $H^{\bullet}(\tilde{Y},\bQ)$ is generated by the set consisting of $e_v:=[E_v]$ and $e_0:=-[\tilde{H}]$. All the multiplications between generators are $0$ except $e_v^2=-e_0^2$ for $v \in V$. The multiplication of any three elements is $0$ since $H^6(\tilde{Y},\bQ)=0$.

By \cite{Fu}-Example 15.4.2, we have
\begin{align*}
c(\tilde{Y})&=(1-e_0)^3\prod_{v\in V}\left( (1+e_v)\left(\frac{1 -e_0- e_v}{1 -e_0}\right)^2\right).
\end{align*}
Thus,
\begin{align*}
c(\Omega_{\tilde{Y}}^1)&=(1+e_0)^3\prod_{v\in V}\left( (1-e_v)\left(\frac{1 +e_0+ e_v}{1 +e_0}\right)^2\right).
\end{align*}
From equation (\ref{log}), we get 
$$c(\Omega_{\tilde{Y}}^1(\log \tilde{Z}))=c(\Omega_{\tilde{Y}}^1)\prod_{v \in V}\frac{1}{1-e_v}\prod_{l \in L}\frac{1}{1-[E_l]}.$$\\
Note that $[E_l]=-\left(d_l e_0+\sum_{v \in V} m_{v,l} e_v \right)$ for $l \in L$ and $v\in V$.

We calculate each factor of the Chern classes above in the ring $H^{\bullet}(\tilde{Y},\bQ)$:
\begin{align*}
\prod_{v \in V}(1-e_v)&=1+\sum_{v \in V}(-e_v),\\
(1+e_0)^3&=1+3 e_0+3 e_0^2,\\
\prod_{v \in V}\left(\frac{1 +e_0+ e_v}{1 +e_0}\right)^2&=\prod_{v \in V}\left(\left(1 +e_0+ e_v\right)\left(1 -e_0+e_0^2\right)\right)^2\\
&=\prod_{v \in V}\left(1 +2e_v+e_v^2 \right) \\
&=1+2\sum_{v \in V} e_v+\sum_{v \in V} e_v^2,
\end{align*}
and 
\begin{align*}
\prod_{l \in L}\frac{1}{1-[E_l]}&=\prod_{l \in L}\frac{1}{1+\left(d_l e_0+\sum_{v \in V} m_{v,l} e_v \right)}\\
&=\prod_{l \in L}\left(1-\left(d_l e_0+\sum_{v \in V} m_{v,l} e_v \right)+\left(d_l e_0+\sum_{v \in V} m_{v,l} e_v \right)^2\right)\\
&=1-\left(\sum_{v \in V} m_{v,red} e_v +d e_0 \right)\\
&+\left(\sum_{v \in V}\frac{m_{v,red}^2+\sum_{l \in L} m_{v,l}^2 }{2} e_v^2+\frac{d^2+\sum_{l \in L}d_l^2}{2}e_0^2\right).
\end{align*}
by $e_0^3=e_v^3=e_0 e_v=e_v e_v'=0$ and $e_v^2=-e_0^2$ for $v,v'\in V$ such that $v\neq v'$.

From these factors we give formulas for the Chern classes:
\begin{align*}
c(\Omega_{\tilde{Y}}^1)&=1+\left(\sum_{v \in V} e_v +3 e_0 \right)+  \left(-\sum_{v \in V}e_v^2 +3e_0^2\right),\\
c(\tilde{Y})&=1-\left(\sum_{v \in V} e_v +3 e_0 \right)+  \left(-\sum_{v \in V}e_v^2 +3 e_0^2\right),
\end{align*}
and
\begin{align*}	
c(\Omega_{\tilde{Y}}^1(\log \tilde{Z}))&=1-\left(\sum_{v \in V} (m_{v,red}-2) e_v +(d-3) e_0\right)\\
	&+\left(\sum_{v \in V}\left(\frac{m_{v,red}^2+\sum_{l \in L} m_{v,l}^2}{2}-2 m_{v,red}+1\right) e_{v}^2  \right. \\
	&\left. +\left(\frac{d_{red}^2+\sum_{l \in L}d_l^2}{2}+3(1-d)\right) e_0^2\right).
\end{align*}

By using Example \ref{wedge} we calculate Chern characters:
 	\begin{align*}
	 	ch(\Omega_{\tilde{Y}}^0(\log \tilde{Z}))&=1,\\
		ch(\Omega_{\tilde{Y}}^1(\log \tilde{Z}))&=2-\left(\sum_{v \in V} (m_{v,red}-2) e_v +(d-3) e_0\right)\\
			&+\left(\sum_{v \in V}\left(\frac{2-\sum_{l \in L }m_{v,l}^2}{2}\right) e_v^2+\left(\frac{3-\sum_{l \in L} d_l^2}{2}\right) e_0^2\right),\\
 	\end{align*}
and
	\begin{align*}	
		ch(\Omega_{\tilde{Y}}^2(\log \tilde{Z}))&=1-\left(\sum_{v \in V} (m_{v,red}-2) e_v +(d-3) e_0 \right)\\
			&+\frac{1}{2}\left(\sum_{v \in V} (m_{v,red}-2)^2 e_v^2 +(d-3)^2 e_0^2\right).\\
	\end{align*}	
Also, we get 
\begin{align*}
ch(\Sf_{\tilde{Y}}(U))&=1+\left(\sum_{v \in V} u_v e_v+u_0 e_0\right)+\frac{1}{2}\left(\sum_{v \in V} u_v^2  e_v^2 + u_0^2 e_0^2\right)
\end{align*}
by formula (\ref{generalch}) and 
\begin{align*}
Td&(\tilde{Y})=1-\frac{1}{2}\left( \sum_{v\in V} e_v+3e_0\right) +e_0^2
\end{align*}
by formula (\ref{td}).

Multiplying all the factors of the right-hand-side of (\ref{Xp}) we get $\X_p(u):=\X(\tilde{Y},\cF_p(U))$:
 	\begin{align*}
		\X_0(u)&=\binom{u_0-1}{2}-\sum_{v \in V} \binom{u_v}{2},\\
		\X_1(u)&=-(u_0-1)(d_{red}-u_0-1) +\sum_{v \in V}u_v(m_{v,red}-u_v-1)\\
		 &-\sum_{l \in L}\binom{d_l}{2} + \sum_{v \in V} \sum_{l \in L}\binom{m_{v,l}}{2},\text{ and}\\
		\X_2(u)&=\binom{d_{red}-u_0-1}{2}-\sum_{v \in V} \binom{m_{v,red}-u_v-1}{2}.
	\end{align*}	
To get $n_{f,k+\frac{j}{d}}$ from (\ref{euler}) we put  
	\begin{align*}
		u&=(d-j) [\tilde{H}]+\sum_{w \in W} \lfloor (d-j)m_w /d\rfloor [E_w]\\
		&=(d-j) [\tilde{H}]+\sum_{v \in V} \lfloor (d-j)m_v /d\rfloor [E_v]+\sum_{l \in L} \lfloor (d-j)m_l /d\rfloor [E_l]\\
		&=(d-j) (-e_0)+\sum_{v \in V} \lfloor (d-j)m_v /d\rfloor e_v-\sum_{l \in L} \lfloor (d-j)m_l /d\rfloor \left(d_l e_0+\sum_{v \in V} m_{v,l} e_v \right)\\
		&=\left( \sum_{l\in L} \lceil j m_l/d \rceil d_l-j\right)e_0+\sum_{v \in V} \left( \sum_{l\in L} \lceil j m_l/d\rceil m_{v,l} - \lceil j m_v/d\rceil \right)e_v.
	\end{align*}	
By denoting  the coefficient of $e_\bullet$ by $u_{\bullet,j}$ we get the formula in Theorem \ref{c3}.


\begin{thebibliography}{EMS}


\bibitem[1]{Bu-S} N. Budur and M.Saito, Jumping coefficients and spectrum of a hyperplane arrangement. Math. Ann. 347 (2010), no.3, 545--579. 17, 23.

\bibitem[2]{Bu-HH}N. Budur, Hodge spectrum of hyperplane arrangements. arXiv:0809.3443(unpublished).

\bibitem[3]{Fu}W. Fulton, Intersection Theory. Springer, Berlin, (1984).

\bibitem[4]{Kul}V. Kulikov, Mixed Hodge structures and singularities. Cambridge Univ. Press, Cambridge, (1998) xxii+186.

\bibitem[5]{Mil}J. Milnor, Singular Points of Complex Hypersurfaces. (AM-61). Princeton University Press, (1968)

\bibitem[6]{St}C.A.M. Peters and J.H.M. Steenbrink, Mixed Hodge Structures. Ergeb. Math. Grenzgeb. 3. Folge Vol. 52. Springer-Verlag Berlin Heidelberg (2008).

\bibitem[7]{YY}Y. Yoon, Spectrum of hyperplane arrangements in four variables. Comm. Algebra {43} (2015), 2585--2600.

\end{thebibliography}
\end{document}